\documentclass[12pt]{amsart}
\usepackage{amsmath,amscd,amssymb,amsfonts}
\def\vers{Aug.~18, 2010, v.2}
\def\nin{\noindent}
\def\bs{\bigskip}
\def\ms{\medskip}
\def\sb{\raise.2ex\h{${\scriptscriptstyle\bullet}$}}
\def\scirc{\,\raise.2ex\hbox{${\scriptstyle\circ}$}\,}
\def\mopl{\hbox{$\bigoplus$}}
\def\msum{\hbox{$\sum$}}
\def\mcup{\hbox{$\bigcup$}}
\def\a{\alpha}
\def\rk{\mathop{\rm rk}}
\def\C{{\mathbf C}}
\def\d{\partial}
\def\D{\Delta}
\def\Db{{\mathbf D}}
\def\Dc{{\mathcal D}}
\def\Dt{\widetilde{\Delta}}
\def\E{{\mathcal E}}
\def\Et{\widetilde{\mathcal E}}
\def\G{\Gamma}
\def\Gt{\widetilde{\Gamma}}
\def\h{\hbox}
\def\Hb{{\mathbf H}}
\def\Hc{{\mathcal H}}
\def\I{{\mathcal I}}
\def\K{{\mathcal K}}
\def\la{\lambda}
\def\Lc{{\mathcal L}}
\def\Lct{\widetilde{\mathcal L}}
\def\M{{\mathcal M}}
\def\N{{\mathbf N}}
\def\Nt{\widetilde{N}}
\def\Oc{{\mathcal O}}
\def\Q{{\mathbf Q}}
\def\q{\quad}

\def\s{\sigma}
\def\t{\widetilde{t}}
\def\Tt{\widetilde{T}}
\def\u{\widehat{u}}
\def\V{{\mathcal V}}
\def\Vt{\widetilde{\mathcal V}}
\def\Z{{\mathbf Z}}
\def\Coim{\hbox{\rm Coim}}
\def\Ker{\hbox{\rm Ker}}
\def\Im{\hbox{\rm Im}}
\def\Ext{{\mathcal E}xt}
\def\Hom{{\mathcal H}om}
\def\Specan{{\mathcal S}pecan}
\def\Gr{\hbox{\rm Gr}}
\def\GGK{{\rm GGK}}
\def\Sch{{\rm Sch}}
\def\1{{\hskip1pt}}
\def\simto{\buildrel\sim\over\longrightarrow}
\def\onto{\mathop{\rlap{$\to$}\hskip2pt\hbox{$\to$}}}
\def\into{\hookrightarrow}

\begin{document}
\title{A variant of N\'eron models over curves}
\author{Morihiko Saito}
\address{RIMS Kyoto University, Kyoto 606-8502 Japan}
\author{Christian Schnell}
\address{Department of Mathematics, Statistics \& Computer Science,
University of Illinois at Chicago
851 South Morgan Street, Chicago, IL 60607}
\date{\vers}
\begin{abstract}
We study a variant of the N\'eron models over curves which is
recently found by the second named author in a more general
situation using the theory of Hodge modules.
We show that its identity component is a certain open subset of
an iterated blow-up along smooth centers of the Zucker extension
of the family of intermediate Jacobians and that the total space
is a complex Lie group over the base curve and is Hausdorff as
a topological space.
In the unipotent monodromy case, the image of the map to the
Clemens extension coincides with the N\'eron model defined by
Green, Griffiths and Kerr.
In the case of families of Abelian varieties over curves, it
coincides with the Clemens extension, and hence with the
classical N\'eron model in the algebraic case (even in the
non-unipotent monodromy case).
\end{abstract}

\maketitle

\section*{Introduction}

\noindent
Let $\Hb$ be a polarizable variation of Hodge structure of weight
$-1$ on a punctured disk $\D^*$.
Let $(\Lc_{\D^*},F)$ be its underlying filtered $\Oc_{\D^*}$-module,
and $(\Lc^{\ge 0},F)$ be its Deligne extension over $\D$ where the
eigenvalues of the residue are contained in $[0,1)$.
Let $\V$ be the vector bundle corresponding to the locally free
sheaf $\Lc^{\ge 0}/F^0\Lc^{\ge 0}$.
Let $\G\subset\V$ denote the subset corresponding to the subsheaf
$j_*\Hb_{\Z}$ where $\Hb_{\Z}$ denotes here the underlying $\Z$-local
system of $\Hb$ and $j:\D^*\into\D$ is the inclusion.
Then the Zucker extension [Zu] of the family of intermediate
Jacobians is defined by
$$J_{\D}^Z(\Hb)=\V/\G.$$
It is a complex analytic Lie group over $\D$, and is the identity
component $J_{\D}^C(\Hb)^0$ of the Clemens extension $J_{\D}^C(\Hb)$,
i.e. the latter is constructed by gluing copies of $J_{\D}^Z(\Hb)$, 
see [Cl], [Sa2].
However, $J_{\D}^Z(\Hb)$ and hence $J_{\D}^C(\Hb)$ are not necessarily
Hausdorff in general.
In the unipotent monodromy case it has been pointed out by M.~Green,
P.~Griffiths and M.~Kerr [GGK2] that it is more natural to consider a
subset $J_{\D}^{\GGK}(\Hb)^0$ of $J_{\D}^Z(\Hb)$ whose fiber over the
origin is the Jacobian of $\Ker\,N\subset H_{\infty}$ where
$H_{\infty}$ is the limit mixed Hodge structure and $N=\log T$ with
$T$ the monodromy.
Indeed, it is a Hausdorff topological space (see [Sa3])
although it is not a complex analytic Lie group over $\D$ in the
usual sense.

Recently the second named author [Sch] found a variant of the
N\'eron model $J_S^{\Sch}(\Hb)$ in a more general situation
where $S$ is a complex manifold which is a partial compactification
of $S^*$ on which the variation of Hodge structure $\Hb$ is defined.
In order to define $J_S^{\Sch}(\Hb)$, consider the polarizable
Hodge module on $S$ naturally extending $\Hb$, and let $(\M,F)$ be
its underlying filtered left $\Dc$-module.
This corresponds by the de Rham functor to the intersection complex
with coefficients in the local system $\Hb_{\C}$, see [Sa1].
Then the identity component $J_S^{\Sch}(\Hb)^0$ is defined by taking
a quotient of the analytic space associated to the symmetric
algebra of $F_0\M$, see [Sch].
In our case where $S=\D$, $\M$ is a $\Dc_{\D}$-submodule of the
union of the Deligne extensions $\Lc^{>-\infty}$, and $F_0\M$ is
a free sheaf on $\D$.
Then $J_{\D}^{\Sch}(\Hb)^0$ is defined by replacing
$\Lc^{\ge 0}/F^0\Lc^{\ge 0}$ in the definition of $J_{\D}^Z(\Hb)$
with the dual of $F_0\M$ which is identified with a free subsheaf
of $\Lc^{\ge 0}/F^0\Lc^{\ge 0}$ using the polarization of $\Hb$.

Set
$$H^{\rm van}_{\infty}=H_{\infty}/H^{\rm inv}_{\infty}\q\h{with}
\q H^{\rm inv}_{\infty}:=\Ker(T-id)\subset H_{\infty}.$$
Here the limit mixed Hodge structure $H_{\infty}$ is defined by
using the base change associated to a cyclic ramified covering of
$\D$, and $H^{\rm inv}_{\infty}$ is a mixed Hodge structure.

\ms\nin
{\bf Theorem~1.} {\it
Let $a=\max\{p\in\N\mid F^pH^{\rm van}_{\infty,\C}\ne 0\}$, and
$d_k=\dim F^kH^{\rm van}_{\infty,\C}$.
There is a sequence of morphisms of complex Lie groups
$\s_k:Y_k\to Y_{k-1}$ over $\D$ for $k=1,\dots,a$ such that
$Y_0=J_{\D}^Z(\Hb)$, $Y_a=J_{\D}^{\Sch}(\Hb)^0$, and $Y_k$ is a
complex manifold defined by $\V_k/\G$ where $\V_k$ are vector
bundles over $\D$ for $k\in[0,a]$.
Moreover, $\s_k$ is obtained by taking the blow-up of $Y_{k-1}$
along a smooth center of codimension $d_k+1$ which is
contained in the fiber $Y_{k,0}$ over $0\in\D$, and by
restricting it to the complement of the strict transform of
$Y_{k-1,0}$.}

\ms\nin
{\bf Theorem~2.} {\it
The image $I_k$ of the composition $\s_1\scirc\cdots\scirc\s_k:Y_k
\to Y_0$ is independent of $k=1,\dots,a$.
Moreover, the unipotent monodromy part of the fiber $(I_k)_0$ of
$I_k$ over $0\in\D$ coincides with the image of
$H^{\rm inv}_{\infty,\C}$ in $J_{\D}^Z(\Hb)_0$.
In particular, $I_k$ for $k>0$ coincides with the identity component
of the N\'eron model $J_{\D}^{\GGK}(\Hb)$ of Green, Griffiths and
Kerr in the unipotent monodromy case.}

\ms
Here the unipotent monodromy part of $(I_k)_0$ means its image by
the surjection from $J_{\D}^Z(\Hb)_0$ to $J_{\D}^Z(\Hb)_0^{\rm unip}
:=H_{\infty,\C,1}/(F^0H_{\infty,\C,1}+H^{\rm inv}_{\infty,\Z})$
where $H_{\infty,\C,1}$ is the unipotent monodromy part of
$H_{\infty,\C}$.
Using an argument on perfect pairings of locally free sheaves,
the proofs are reduced to the calculation of $F_0\Gr^{\a}_V\M$
where $V$ is the filtration of Kashiwara and Malgrange, and
this is calculable in terms of the limit mixed Hodge structure by
the theory of Hodge modules ([Sa1], 3.2).
Using Theorems~1 and 2 together with [Sa2], [Sa3], we get moreover

\ms\nin
{\bf Theorem~3.} {\it
With the above notation, the $Y_k$ for $k\in[0,a]$ and hence
$Y_a=J_{\D}^{\Sch}(\Hb)$ are complex Lie groups over $\D$.
Moreover these are Hausdorff topological spaces if $k>0$.}

\ms
For $k=0$, $Y_0=J_{\D}^Z(\Hb)$ is not necessarily Hausdorff,
although it is Hausdorff on a neighborhood of the image of $\s_1$,
see [Sa2], [Sa3].
Note that Hausdorff property is not included in the definition
of complex manifold in loc.~cit.\ and also in this paper.
From Theorem~1 we can deduce

\ms\nin
{\bf Corollary.} {\it
In case of families of Abelian varieties, $J_{\D}^{\Sch}(\Hb)$
coincides with the Clemens extension $J_{\D}^C(\Hb)$, and hence with
the classical N\'eron model in the algebraic case.}

\ms
Indeed, the first assertion follows from Theorem~1 since
$a=0$ in this case.
The last assertion follows from [Sa2], 4.5.
Here we do not assume the monodromy unipotent.

\ms
To illustrate Theorem~3, we describe what happens for variations of
Hodge structure on $\D^*$ of ``mirror quintic type''
up to a Tate twist.
This means that the rank of the local system is 4 and the nonzero
Hodge numbers of the general fibers are given by $h^{p,q} = 1$ for 
$p=-2,-1,0,1$ with $q=-p-1$.
We assume that the monodromy around the origin is unipotent.
In this case, $\V$ is a vector bundle of rank 2, and the central fiber
$\V_0$ is $H_{\infty,\C}/F^0 H_{\infty,\C}$.

In the notation of [GGK1], there are three types of degenerations.
When $N^2 = 0$ and $\rk N = 1$ (type ${\rm II}_1$), we have
$F^0 H_{\infty,\C}^{\rm van} = 0$, i.e. $\Ker\,N$ surjects onto
$\V_0/\G_0$, and so $J_{\D}^{\Sch}(\Hb)^0 = J_{\D}^Z(\Hb)$.
When $N^2 = 0$ and $\rk N = 2$ (type ${\rm II}_2$), or when
$N^3 \neq 0$ (type ${\rm I}$), the subspace
$F^0 H_{\infty,\C}^{\rm van}$ is one-dimensional, i.e. the image of
$\Ker\,N$ in $J_{\D}^Z(\Hb)_0=\V_0/\G_0$ has codimension 2 in
$J_{\D}^Z(\Hb)$.
So $J_{\D}^{\Sch}(\Hb)^0$ is obtained from the Zucker extension
$J_{\D}^Z(\Hb)$ by blowing up it along the image of $\Ker\,N$ in
$\V_0/\G_0$ and deleting the strict transform of the central fiber
$\V_0/\G_0$ in this case.

In the non-unipotent monodromy case, the situation is similar.
Tensoring the above examples with the local system of rank 1 with
monodromy $-1$, we get examples with non-unipotent monodromy.
We get $J_{\D}^{\Sch}(\Hb)^0$ in the same way as above except that
the center of the blow-up in the last two cases is not the image of
$\Ker\,N$.

\ms
Finally we note a few remarks on $J_S^{\Sch}(\Hb)$ and other N\'eron
models in the higher dimensional case:

For families of Abelian varieties defined outside a divisor with
normal crossings, $F_0\M$ is a free subsheaf of a Deligne
extension whose quotient is also free, and
$J_S^{\Sch}(\Hb)^0$ coincides with the generalized Zucker extension.
A generalization of the N\'eron model in this case has been constructed 
by A.~Young [Yo] assuming the local monodromies unipotent and using
a different method.

In general $F_0\M$ is not necessarily free or reflexive
even in the normal crossing case, and $\sigma^*F_0\M$ may
have torsion for a blow-up $\sigma$.
If we replace $F_0\M$ with its reflexive hull, i.e.\ the double
dual $(F_0\M)^{\vee\vee}$, then the latter is reflexive and the
morphism $(F_0\M)^{\vee}\to(F_0\M)^{\vee\vee\vee}$ is an isomorphism,
see e.g.\ [OSS].
In this way, it may be possible to extend some of the arguments using
the pairings in this paper at least to the normal crossing case.

It is shown in [Sch] that there is a natural surjection (a kind of
`blow-down') from $J_S^{\Sch}(\Hb)$ onto the generalized N\'eron
model $J_S^{\rm BPS}(\Hb)$ defined in [BPS]; for families of Abelian
varieties over curves, this was noted in [BPS], Remark~2.7(i).

\ms
We would like to thank the referee for useful comments.
The first named author is partially supported by Kakenhi 21540037.

Here is a brief outline of the paper:
In Section 1 we review some facts related to vector bundles and
locally free sheaves over curves including Deligne extensions and
$V$-filtrations.
In Section 2 we prove the main theorems using these.
In Section 3 we give some remarks for the case where the base
space $S$ is not a curve.

\bs\bs
\section{Preliminaries}

\nin
{\bf 1.1.~Vector bundle case.}
We first consider the vector bundle case (before dividing out by
$\G$ in the introduction).
In general, let $\E$ be a free sheaf of rank $r$ on $\D$, and
let $\V(\E)$ be the corresponding vector bundle over $\D$ such that
$$\E=\Oc_{\D}(\V(\E)),\q\V(\E)=\Specan_{\D}
({\rm Sym}_{\Oc}^{\sb}\E^{\vee}),\leqno(1.1.1)$$
where $\Oc_{\D}(\V(\E))$ is the sheaf of local sections of $\V(\E)$,
$\E^{\vee}:=\Hom_{\Oc}(\E,\Oc_{\D})$ is the dual free sheaf of $\E$,
and ${\rm Sym}_{\Oc}^{\sb}\E^{\vee}$ is the symmetric algebra.

Let $\E'$ be a free subsheaf of $\E$ such that $\E/\E'$ is
$t$-torsion, where $t$ is the coordinate of $\D$.
Then there are nonnegative integers $a_i$ such that
$$\E/\E'=\mopl_{i=1}^r\bigl({\bf C}[t]/t^{a_i}{\bf C}[t]\bigr).$$
More precisely, by the theory of modules over principal ideal
domains, there are integers $a_i$ and bases $v_1,\dots,v_r$ and
$v'_1,\dots,v'_r$ of $\E$ and $\E'$ respectively (shrinking $\D$
if necessary) such that
$$v'_i=t^{a_i}v_i\q\h{and}\q a_i\ge a_{i+1}\ge 0.\leqno(1.1.2)$$
Let $x_1,\dots,x_r$ and $y_1,\dots,y_r$ be respectively the
coordinates of the vector bundles $\V(\E)$, $\V(\E')$ associated to
the bases $v_1,\dots,v_r$ and $v'_1,\dots,v'_r$, i.e.\
$\langle x_i,v_j\rangle=\delta_{i,j}$ and similarly for $y_i, v'_j$,
where the $x_i$ and $y_i$ are identified with sections
of $\E^{\vee}$ and $\E'{}^{\vee}$.
Then
$$x_i=t^{a_i}y_i.\leqno(1.1.3)$$
Set $a=\max\{a_i\}$, and $J_k=\{i\mid a_i\ge k\}\,\,(k=0,\dots,a)$.
Then $J_k=\{1,\dots,m_k\}$ with $m_k\ge m_{k+1}$. Define
$$\E_k=\E'+t^k\E\,\,(k=0,\dots,a).\leqno(1.1.4)$$
Then $\E_0=\E$, $\E_a=\E'$, and $\E_k$ is the free subsheaf of
$\E$ generated by
$$t^{c_{k,i}}v_i\q\h{with}\q c_{k,i}=\min(a_i,k)\,\,(i\in[1,r]).$$
Note that $c_{k,i}\le k$, and $c_{k,i}=k$ if and only if $i\le m_k$.
Let $x_1^{(k)},\dots,x_r^{(k)}$ be the coordinates of the vector
bundle $\V(\E_k)$ corresponding to the basis
$t^{c_{k,1}}v_1,\dots,t^{c_{k,r}}v_r$.
Then for $k\in[1,a]$
$$x_i^{(k-1)}=\begin{cases}tx_i^{(k)}&\text{if $\,i\le m_k$,}\\
x_i^{(k)}&\text{if $\,i>m_k$.}\end{cases}\leqno(1.1.5)$$
Thus $\V(\E_k)$ is an open subset of the blow-up of $\V(\E_{k-1})$
along a center of dimension $r-m_k$.
More precisely, $\V(\E_k)$ is the complement of the proper
transform of the fiber $\V(\E_{k-1})_0$ over $0\in\D$.

\ms\nin
{\bf 1.2.~Duality.}
Let $\E$ be a free sheaf on $\D$, and $\K:=\E[t^{-1}]$ be the
localization by the coordinate $t$ of $\D$.
In this case we say that $\E$ is a lattice of $\K$.
We will fix $\K$ and consider lattices $\E$ of $\K$.

In this subsection, we denote the dual free sheaf by
$$\Db(\E):=\Hom_{\Oc}(\E,\Oc_{\D}).$$
There is a perfect pairing
$$\langle *,*\rangle:\Db(\E)\otimes_{\Oc}\E\to\Oc_{\D},$$
corresponding to the canonical isomorphism
$$\Db(\E)\simto\Hom_{\Oc}(\E,\Oc_{\D}).$$
The above pairing is extended to
$$\langle *,*\rangle:\Db(\K)\otimes_{\Oc}\K\to\Oc_{\D}[t^{-1}],$$
where
$$\Db(\K):=\Hom_{\Oc[t^{-1}]}(\K,\Oc_{\D}[t^{-1}]).$$
(Here $\K$ is not a $\Dc$-module in general, and this is quite
different from the dual as a $\Dc$-module in [Sa1] in case
$\K$ has a structure of a $\Dc$-module.)
We have
$$\Db(\Db(\K))=\K,\q\Db(\Db(\E))=\E,\leqno(1.2.1)$$
$$\Db(\E)=\{\xi\in\Db(\K)\mid\langle\xi,v\rangle\in\Oc_{\D}\,\,\,
\h{for any}\,\,v\in\E\}.\leqno(1.2.2)$$
Note also that for any basis $v_1\dots,v_r$ of $\E$, there is a
unique dual basis $v_1^*,\dots,v_r^*$ of $\Db(\E)$ such that
$$\langle v_i^*,v_j\rangle=\delta_{i,j}.$$

Using (1.2.1--2), we have for lattices $\E,\E'$ of $\K$
$$\Db(\E)\cap\Db(\E')=\Db(\E+\E'),\q
\Db(\E)+\Db(\E')=\Db(\E\cap\E').\leqno(1.2.3)$$

For a $t$-torsion $\Oc_{\D}$-module $\E''$, define
$$\Db^1(\E''):=\Ext_{\Oc}^1(\E'',\Oc_{\D}).$$
This is defined by taking a free resolution
$0\to E_1\to E_0\to\E''\to 0$.
Note that
$$\Ext_{\Oc}^i(\E'',\Oc_{\D})=0\q\h{for}\,\,\,i\ne 1.$$

For a short exact sequence
$$0\to\E'\to\E\to\E''\to 0,$$
such that $\E',\E$ are lattices (and hence $\E''$ is
$t$-torsion), we have the dual exact sequence
$$0\to\Db(\E)\to\Db(\E')\to\Db^1(\E'')\to 0.\leqno(1.2.4)$$
Moreover, if $\E''$ is annihilated by $t$, then we have an
isomorphism as $\C$-vector spaces
$$\Db^1(\E'')=\Db_{\C}(\E''):=\Hom_{\C}(\E'',\C),\leqno(1.2.5)$$
using a free resolution $0\to\Oc_{\D}\otimes\E''\buildrel{t}\over\to
\Oc_{\D}\otimes\E''\to\E''\to 0$.

\ms\nin
{\bf 1.3.~Relation between $F$ and $V$.}
Let $(\M,F)$ be a filtered left $\Dc_{\D}$-module underlying a
polarizable Hodge module with strict support $\D$.
The condition on strict support is equivalent to the condition that
$\M$ has no nontrivial sub nor quotient module supported on points,
and is further equivalent to the condition that it corresponds by
the de Rham functor to an intersection complex (see [BBD]) with
local system coefficients, see [Sa1], 5.1.3.
Let $\M[t^{-1}]$ be the localization of $\M$ by the coordinate $t$
of $\D$.

Shrinking $\D$ if necessary we may assume that the restriction
$\Lc_{\D^*}:=\M|_{\D^*}$ is locally free over $\Oc_{\D^*}$, i.e.
$\Lc_{\D^*}$ corresponds to a local system on $\D^*$.
Let $\Lc^{\ge\a}$ (resp. $\Lc^{>\a}$) denote the Deligne extension
of $\Lc_{\D^*}$ such that the eigenvalues of the residue of the
connection are contained in $[\a,\a+1)$ (resp. $(\a,\a+1]$), see [D1].
These are identified with the $V$-filtration of Kashiwara and
Malgrange on $\M[t^{-1}]$ indexed by ${\bf Q}$.
Here we use left $\Dc$-modules, and assume that the action of
$t\d_t-\a$ is nilpotent on $\Gr_V^{\a}\M$.
Note that $\M[t^{-1}]$ is identified with the union of Deligne
extensions $\Lc^{>-\infty}$, and there are canonical isomorphisms
$$\Lc^{\ge\a}=V^{\a}\Lc^{>-\infty},
\q\Lc^{>\a}=V^{>\a}\Lc^{>-\infty}.$$
It is also well-known that $F^p\Lc_{\D^*}$ are extended as locally
free subsheaves $F^p\Lc^{\ge\a}$ of $\Lc^{\ge\a}$ such that the
$\Gr^F_p\Lc^{\ge\a}$ are locally free for any $\a\in\Q$.
We have similar assertions with $\Lc^{\ge\a}$ replaced by $\Lc^{>\a}$.

There are canonical inclusions
$$\Lc^{>-1}\into\M\into\Lc^{>-\infty},$$
and moreover we have by the definition of Hodge modules
([Sa1], 3.2)
$$F_p\M=\msum_{i\ge 0}\,\d_t^i(F^{i-p}\Lc^{>-1}),
\leqno(1.3.1)$$
and $\M=\Dc_{\D}\Lc^{>-1}$ for $p\to\infty$, where
$\d_t:=\d/\d t$ and $F_p=F^{-p}$.
Note that (1.3.1) implies
$$F_p\M\cap\Lc^{>-1}=F^{-p}\Lc^{>-1},\leqno(1.3.2)$$
and assuming (1.3.2), condition (1.3.1) is equivalent to
the surjections (see [Sa1], 3.2.2)
$$\d_t:F_p\Gr_V^{\a}\M\onto F_{p+1}\Gr_V^{\a-1}\M
\q\h{for any}\,\,\,p\in\Z,\,\a\le 0.\leqno(1.3.3)$$
Note also that the morphism (1.3.3) for $\a<0$ is injective (and
hence bijective) since it is an isomorphism if $F_p$ and $F_{p+1}$
are omitted.
So we get the isomorphisms
$$\d_t:F_p\Gr_V^{\a}\M\simto F_{p+1}\Gr_V^{\a-1}\M
\q\h{for any}\,\,\,p\in\Z,\,\a<0.\leqno(1.3.4)$$
Then (1.3.3) is equivalent to (1.3.4) together with
$$\d_t:\Gr_V^0(\M,F)\to\Gr_V^{-1}(\M,F[-1])\,\,\,\h{is strictly
surjective,}
\leqno(1.3.5)$$
where $(F[k])_p=F_{p-k}$.

\ms\nin
{\bf 1.4.~Relation with $(H_{\infty,\C},F)$.}
By definition, $H_{\infty,\C}$ is identified with
$$\Lc^{\ge 0}(0)=\Lc^{\ge 0}/\Lc^{\ge 1}.$$
However, this is not compatible with the Hodge filtration $F$ in
general.
We have to use the following isomorphism (which is a special case
of [Sa1], 3.4.12)
$$H_{\infty,\C}=\mopl_{\a\in[0,1)}\,\Gr_V^{-\a}\Lc^{>-1}=
\mopl_{\a\in[0,1)}\,\Gr_V^{-\a}\M,
\leqno(1.4.1)$$
where the last isomorphism follows from (1.3.2).
Set
$$H_{\infty,\C,\la}:=\Ker(T_s-\la)\subset H_{\infty,\C},$$
where $T_s$ is the semisimple part of the monodromy $T$.
Then (1.4.1) induces
$$(H_{\infty,\C,\la},F)=\Gr_V^{-\a}(\M,F)\q\h{for}\,\,\,
\a\in[0,1),\,\,\la=e^{2\pi i\a},
\leqno(1.4.2)$$
using the unipotent base change as in (1.5) below, where $F_p=F^{-p}$.

Combined with (1.3.3) (or (1.3.4--5)), (1.4.2) induces the
isomorphisms
$$(H^{\rm van}_{\infty,\C,\la},F)\simto\Gr_V^{-\a+j}(\M,F[-j])
\q\h{for}\,\,\,\a\in[0,1),\,\,\la=e^{2\pi i\a},\,\,j\in\Z_{>0}.
\leqno(1.4.3)$$
Indeed, this is clear for $\a\in(0,1)$, and the case $\a=0$ is
reduced to the case $j=1$ by (1.3.4).
In the last case, $-N/2\pi i$ is identified with the composition of
$$\d_t:\Gr_V^0\M\to\Gr_V^{-1}\M\q\h{and}\q t:\Gr_V^{-1}\M\to\Gr_V^0\M,
\leqno(1.4.4)$$
and the kernel of the first morphism is identified by (1.4.2) with
$\Ker\,N\subset H_{\infty,\C,1}$ since the last morphism of (1.4.4)
is injective, see [Sa1], 5.1.3.
Then (1.4.3) follows in this case from the strict surjectivity of
(1.3.5) together with (1.4.2).

\ms\nin
{\bf 1.5.~Deligne extension and the unipotent base change.}
Let $L$ be a local system on $\D^*$ with a quasi-unipotent monodromy,
and $\Lc^{\ge 0}$ be the Deligne extension of $L$ such that the
eigenvalues of the residue of the connection are contained in $[0,1)$.
Let $T=T_uT_s$ be the Jordan decomposition of the monodromy $T$
where $T_s$ and $T_u$ are respectively the semisimple and
unipotent part of $T$.
Set $N=\log T_u$.
For a multivalued section $u$ of $L$, let $u=\msum_j\,u_j$ be
the decomposition such that $T_su_j=\exp(-2\pi i\a_j)u_j$ with
$\a_j\in[0,1)$.
Then we have a corresponding section $\u$ of the Deligne extension
$\Lc^{\ge 0}$ defined by
$$\u=\msum_j\,\exp\Bigl(-\frac{\log t}{2\pi i}N\Bigr)t^{\a_j}u_j.
\leqno(1.5.1)$$

Let $\pi:\Dt\to\D$ be a unipotent base change.
By definition it is an $m$-fold ramified covering of open disks such
that $\pi^*t=\t^m$ and the monodromy $\Tt$ on $\Dt$ is unipotent,
i.e.\ $T_s^m=id$, where $t$ and $\t$ are the coordinates of $\D$ and
$\Dt$ respectively.
Set $\Nt=\log\Tt$.
Since $\Nt$ corresponds to $mN=m\log T_u$, we get
$$\pi^*\u=\msum_j\,\exp\Bigl(-\frac{\log \t}{2\pi i}\Nt\Bigr)
\t^{\,m\a_j}u_j,\leqno(1.5.2)$$
where $m\a_j\in\N$ by hypothesis.
This implies that $\pi^*\Lc^{\ge 0}$ is naturally identified with a
subsheaf of the Deligne extension $\Lct^{\ge 0}$, and the
$V$-filtration on
$\Lc^{\ge 0}$ is induced by the $\t$-adic filtration on $\Lct^{\ge 0}$.

If $L$ underlies a polarizable variation of Hodge structure of
weight $-1$, then we can define as in the introduction
$$\E=\Lc^{\ge 0}/F^0\Lc^{\ge 0},\q\G\subset\V:=\V(\E).$$
We can repeat the same for the pullback to $\Dt$, and get
$$\Et=\Lct^{\ge 0}/F^0\Lct^{\ge 0},\q\Gt\subset\Vt:=\V(\Et).$$
If $\pi^*$ denote also the base change by the morphism $\pi$, then
we have canonical morphisms
$$\pi^*\G\into\Gt,\q \pi^*\V\to\Vt,$$
since $\pi^*\E$ is a subsheaf of $\Et$ and $\pi^*\V$ is associated
to $\pi^*\E$ as in (1.1.1).
Moreover, $\G$ and $\V$ are obtained by taking the quotient of
$\pi^*\G$ and $\pi^*\V$ by the action of the covering transformation
group $G$.
So the conditions (2.3.1--2) for $\G\subset\V$ are reduced to those
for $\Gt\subset\Vt$.
Moreover, (2.3.2) is satisfied for $X$ if it is satisfied for
the pullback of $X\cap\Ker(T_s-id)$.
Indeed, the non-unipotent monodromy part causes no problem by
(1.5.2) since $m\a_j\ge 1$ if $\a_j\ne 0$.

\section{Proof of the main theorems}

\nin
{\bf 2.1.~Construction.}
Let $\Hb$ be a polarizable variation of Hodge structure of weight
$-1$ on $\D^*$, and $\Lc^{\ge\a},\Lc^{>\a}$ be the Deligne extensions
as in (1.3).
Let $(\M,F)$ be the filtered left $\Dc_{\D}$-module underlying
the polarizable Hodge module extending $\Hb$ over $\D$.
This is determined by $\Lc^{>-1}$ as in (1.3.1).
In this paper we use left $\Dc$-modules whereas right $\Dc$-modules
are used in [Sa1].
The transformation between left and right $\Dc_X$-modules on a
complex manifold $X$ is given by
$$(\M,F)\mapsto(\omega_X,F)\otimes_{\Oc}(\M,F)$$ for filtered left
$\Dc_X$-modules $(\M,F)$, where $\omega_X=\Omega_X^{\dim X}$ is
the dualizing sheaf of $X$ and the Hodge filtration $F$ on
$\omega_X$ is defined by $\Gr_F^p\omega_X=0$ for $p\ne -\dim X$.

Let $\Db(F_0\M)$ be the dual sheaf as in (1.2), and
$\V\bigl(\Db(F_0\M)\bigr)$ be the associated vector bundle over
$\D$ as in (1.1).
Here $F_0\M\subset\Lc^{>-\infty}$ is torsion-free, and is hence a
free $\Oc_{\D}$-module.

Using a polarization of the variation of Hodge structure $\Hb$,
we have perfect pairings for $\a\in\Q$
$$\aligned \Lc^{\ge\a}&\otimes_{\Oc}\Lc^{>-\a-1}\to\Oc_{\D},\\
\bigl(\Lc^{\ge\a}/F^0\Lc^{\ge\a}\bigr)&\otimes_{\Oc}F^0\Lc^{>-\a-1}
\to\Oc_{\D}.\endaligned$$
Define
$$\E^{\ge\a}:=\Lc^{\ge\a}/F^0\Lc^{\ge\a},\q\E:=\E^{\ge 0}.
\leqno(2.1.1)$$
In the notation of (1.2), the above perfect pairings imply the
identifications
$$\Db(\Lc^{>-\a-1})=\Lc^{\ge\a},\q
\Db(F^0\Lc^{>-\a-1})=\E^{\ge\a},$$
and similarly with $>-\a-1$ and $\ge\a$ replaced respectively by
$\ge -\a-1$ and $>\a$.

Since (1.3.1) implies
$$F^0\Lc^{>-1}\subset F_0\M,$$
we have by (1.2.4) the inclusion
$$\E' = \Db(F_0\M)\subset\E.$$

We define the identity component $J_{\D}^{\Sch}(\Hb)^0$ by
$$J_{\D}^{\Sch}(\Hb)^0=\V\bigl(\Db(F_0\M)\bigr)/\G,$$
where $\G$ is the subspace of $\V\bigl(\Db(F_0\M)\bigr)$
corresponding to the subsheaf
$$j_*\Hb_{\Z}\into\Db(F_0\M)\subset\E:=\Lc^{\ge 0}/F^0\Lc^{\ge 0}.$$
Here we have the inclusion over $\D^*$ since $\Hb$ has weight $-1$.
For the inclusion over $0\in\D$, we have to show
$$\langle u,v\rangle\in\Oc_{\D}\q\h{for any}\,\,\,u\in j_*\Hb_{\Z},\,
v\in F_0\M.\leqno(2.1.2)$$
Since we have the injection $j_*\Hb_{\Z}\into\E$ and $\d_tu$
vanishes, (2.1.2) follows from (1.3.1) using
$$\d_t\langle u,v\rangle=\langle\d_tu,v\rangle+
\langle u,\d_tv\rangle\q\h{for any }\,\,u,v\in\Lc^{>-\infty},
\leqno(2.1.3)$$
where $\Lc^{>-\infty}=\M[\frac{1}{t}]=\Lc^{>-1}[\frac{1}{t}]$.
(It is shown in [Sch] that (2.1.2) holds for any $u\in j_*\Hb_{\Z}$
and $v\in\M$ in a more general situation.)

Set
$$G=H^1(\D^*,\Hb_{\Z})_{\rm tor}.$$
For any $g\in G$, let $\nu_g$ be an admissible normal function
whose cohomology class is $g$, see [Sa2].
We can then define $J_{\D}^{\Sch}(\Hb)$ in the same way as in
loc.~cit., i.e.
$$J_{\D}^{\Sch}(\Hb)=\mcup_{g\in G}\,J_{\D}^{\Sch}(\Hb)^g
\q\h{with}\q J_{\D}^{\Sch}(\Hb)^g:=\nu_g+J_{\D}^{\Sch}(\Hb)^0.
\leqno(2.1.4)$$
By Proposition~(2.2) below, the $g$-component $J_{\D}^{\Sch}(\Hb)^g$
is independent of the choice of $\nu_g$.
Since $G$ is torsion, we may assume here
$$\nu_g\in\bigl(\h{$\frac{1}{m}$}\G|_{\D^*}\bigr)\big/
\G|_{\D^*}\q\h{for some}\,\,\,m\in\N.\leqno(2.1.5)$$

\ms
The following proposition is a refinement of [EZ], and is proved in
a more general situation in [Sch] using the theory of duality of
mixed Hodge modules [Sa1].
We give here a proof in the curve case using (1.2--4).

\ms\nin
{\bf Proposition~2.2.} {\it
Let $\nu$ be any admissible normal function whose cohomology class
vanishes.
Then it extends to a section of $J_{\D}^{\Sch}(\Hb)^0$.}

\ms\nin
{\it Proof.}
Corresponding to an admissible normal function $\nu$,
we have a short exact sequence of admissible variations of
mixed Hodge structures
$$0\to\Hb\to\Hb^e\to\Z_{\D^*}\to 0,$$
and the cohomology class of $\nu$ is defined by the extension
class of the underlying $\Z$-local systems.
This short exact sequence is easily extended to a short exact
sequence of mixed Hodge modules since $H$ has weight $-1$
(and the intermediate extension of perverse sheaves corresponds
to the minimal extension of $\Dc$-modules).
In particular, we have a short exact sequence of underlying
filtered $\Dc_{\D}$-modules
$$0\to(\M,F)\to(\M^e,F)\buildrel{\rho}\over\to
(\Oc_{\D},F)\to 0.$$
Note that $\rho$ is bistrictly compatible with $F,V$
where $V$ is the filtration of Kashiwara and Malgrange indexed by
$\Q$.
Indeed, $\Gr_V^{\a}\rho$ underlies a morphism of mixed Hodge
modules and hence is strictly surjective, see [Sa1], 3.3.3.
So we have a splitting $\s_F$ which preserves $F$ and $V$
(shrinking $\D$ if necessary).
On the other hand, we have a splitting $\s_{\Z}$ of $\rho$
defined over $\Z$ since the cohomology class of $\nu$ vanishes.
It preserves the filtration $V$ since it is a morphism of
$\Dc$-modules.
Thus the normal function $\nu$ is represented by
$$\nu':=\s_F(1)-\s_{\Z}(1)\in\Lc^{\ge 0}.$$
We have to show $\nu'\in\Db(F_0\M)$, i.e.
$$\langle\nu',v\rangle\in\Oc_{\D}\q\h{for any}\,\,\,v\in F_0\M.$$
By the Griffiths transversality we have
$$\d_t^i\nu'|_{\D^*}\in F^{-i}\Lc_{\D^*}\q\h{for}\,\,\,i>0.$$
Since $\Hb$ has weight $-1$, we have by the definition of
polarization
$$\langle F^i\Lc_{\D^*},F^j\Lc_{\D^*}\rangle=0\q\h{if}\,\,\,
i+j\ge 0.$$
So the assertion follows from (1.3.1) using (2.1.3).

\ms\nin
{\bf Remarks~2.3.} (i)
With the notation of the introduction, the construction in (1.1)
is compatible with taking the quotient by $\G$ if $\G\subset\V(\E')$
and the following condition is satisfied.
$$\h{There is an open neighborhood $U$ of $0\in\V_0$ in $\V$
such that $U\cap\G\subset 0_{\D}$.}\leqno(2.3.1)$$
Here $0_{\D}$ denotes the zero section.
Condition~(2.3.1) is preserved by blowing-ups as in
Theorem~1.
Note that (2.3.1) is equivalent to the condition that, for any point
$p$ of $\V_0$, there exists a neighborhood $U_p$ of $p$ in $\V$ such
that the map $U_p\to\V/\G$ is injective.
Indeed, if (2.3.1) is satisfied, then take $U_p$ such that
$U_p-U_p\subset U$.
If the latter condition is satisfied, then $U$ can be taken to be the
intersection of $U_0$ with the pull-back of $U_0\cap 0_{\D}$.

So (2.3.1) is equivalent to the condition that the quotient has the
induced structure of a complex manifold (or a complex analytic Lie
group over the base space in this case) although the Hausdorff
property is unclear.
(In this paper, a complex manifold means a ringed space which
is locally isomorphic to $(\D^n,\Oc_{\D^n}$), and the Hausdorff
property is not assumed.)
A similar argument is noted in [Sa2], Remark after 3.4.

\ms
(ii) With the above notation, let $X$ be a vector subspace of $\V_0$.
Then the Hausdorff property for any distinct two points of
$X/(\G_0\cap X)$ in $\V/\G$ is equivalent to the following.
$$\h{For any $p\in X\setminus\G_0$, there is an neighborhood
$U_p$ in $\V$ with $U_p\cap\G=\emptyset$.}\leqno(2.3.2)$$
Indeed, the Hausdorff property for the images of $p_1,p_2\in X$ in
$\V/\G$ is equivalent to the existence of neighborhoods $U_{p_i}$
of $p_i$ in $\V$ $(i=1,2)$ such that
$(U_{p_1}-U_{p_2})\cap\G=\emptyset$.
If (2.3.2) is satisfied, then take neighborhoods $U_{p_i}$
such that $U_{p_1}-U_{p_2}\subset U_p$ where $p=p_1-p_2$.
If the latter condition is satisfied, then apply this to $p$ and $0$
so that $U_p\cap\G\subset(U_p-U_0)\cap\G=\emptyset$ where we may
replace $U_p$ by its intersection with the pull-back of $U_0\cap 0_{\D}$.

The combination of (2.3.1) and (2.3.2) is thus equivalent to
the condition that $\V/\G$ has the induced structure of a Hausdorff
complex manifold on a neighborhood of $X/X\cap\G_0$,
and is equivalent to the following single condition:
$$\h{For any $p\in X$, there is an neighborhood
$U_p$ in $\V$ with $U_p\cap\G\subset\G(p)$.}\leqno(2.3.3)$$
Here $\G(p)$ is the section of $\G$ passing through $p\in\V_0$
if $p\in\G$, and is empty otherwise.

Note that (2.3.3) is also stable by blowing-ups as in (1.1).
More precisely, (2.3.3) is satisfied by $\V(\E_k)$ in (1.1) if
(2.3.3) is satisfied on a neighborhood of the image of $\V(\E_k)$
in $\V(\E_{k-1})$.

Using the above equivalence together with the finiteness of the
cohomology classes of admissible normal functions in the curve case,
we can reduce the proof of the property that $J_{\D}^{\Sch}(\Hb)$ is
a Hausdorff complex manifold to that for $J_{\D}^{\Sch}(\Hb)^0$.
Indeed, the assertion is clear for the property that
$J_{\D}^{\Sch}(\Hb)$ is a complex manifold by definition.
As for the Hausdorff property note that the normal function $\nu_g$
in (2.1.4) is represented by a multivalued section $u$ of
$\frac{1}{m}\G|_{\D^*}$ as in (2.1.5).
Hence any point $q$ of the $g$-component $\nu_g+J_{\D}^{\Sch}(\Hb)^0$
over $0\in\D$ is formally represented by $u+p$ with $p\in\V_0$.
This means that if we consider a neighborhood $U_p$ of $p$ in $\D$,
then the restriction over $\D^*$ of the corresponding neighborhood
of $q$ is represented by $u+U_p|_{\D^*}$.
For the proof of the Hausdorff property, it is then sufficient to
show the following:

\ms
For any $p_1,p_2\in\V_0$ and for any multivalued sections $u_1,u_2$
of $\frac{1}{m}\G|_{\D^*}$ such that
$$u_1-u_2\notin\G|_{\D^*}\q\h{or}\q p_1-p_2\notin\G,$$
there are respectively open neighborhoods $U_1,U_2$ of $p_1,p_2$ in
$\V$ such that
$$(U_1-U_2)\cap(u_2-u_1+\G|_{\D^*})\subset\G(p_1-p_2).$$
(Note that the two points represented by $u_1+p_1$ and $u_2+p_2$
are in the same $g$-component if and only if $u_1-u_2\in\G|_{\D^*}$.)
Here we may replace
$$u_1-u_2+\G|_{\D^*}\q\h{with}\q\h{$\frac{1}{m}$}\G,$$
since $u_1-u_2\subset\frac{1}{m}\G|_{\D^*}$.
We may then replace further
$$\G(p_1-p_2)\q\h{with}\q\bigl(\h{$\frac{1}{m}$}\G\bigr)(p_1-p_2).$$
So the proof of the Hausdorff property for $J_{\D}^{\Sch}(\Hb)$ is
reduced to the case $m=1$ by replacing $\G$ with $\frac{1}{m}\G$,
and then follows from (2.3.3) if $J_{\D}^{\Sch}(\Hb)^0$ has the
induced structure of a Hausdorff complex manifold.

\ms\nin
{\bf 2.4.~Proof of Theorem~1.}
Let $V$ denote also the quotient filtration on
$$\E:=\Lc^{\ge 0}/F^0\Lc^{\ge 0}.$$
By [Sa1], (3.2.1.2), we have in the notation of (2.1.1)
$$V^{\a}\E=\E^{\ge\a},\q V^{>\a}\E=\E^{>\a}
\q(\a\ge 0).$$
This implies
$$\Gr_V^{\a}(\E/\Db(F_0\M))=
\frac{\bigl(\Db(F_0\M)+\E^{\ge\a}\bigr)/\Db(F_0\M)}
{\bigl(\Db(F_0\M)+\E^{>\a}\bigr)/\Db(F_0\M)},$$
Considering the definition (1.1.4), it is then sufficient to show
$$d_{j+1}=\msum_{0\le\a<1}\dim\Gr_V^{\a+j}
\bigl(\E/\Db(F_0\M)\bigr)\q(j\in\N).\leqno(2.4.1)$$
Indeed, we have
$$tV^{\a}\E=V^{\a+1}\E\,\,\,(\a\ge 0),$$
and hence $V^{\a}\E$ is a refinement of the $t$-adic filtration
$V^i\E=t^i\E\,(i\in\N)$.

By (1.2.3) we have
$$\aligned &\Db(F_0\M)+\E^{\ge\a}=
\Db(F_0\M\cap F^0\Lc^{>-\a-1})=\Db(F_0V^{>-\a-1}\M),\\
&\Db(F_0\M)+\E^{>\a}=
\Db(F_0\M\cap F^0\Lc^{\ge-\a-1})=\Db(F_0V^{-\a-1}\M).\endaligned$$
Combining these with (1.2.4--5), we get thus
$$\Db_{\C}(F_0\Gr_V^{-\a-1}\M)=\Gr_V^{\a}(\E/\Db(F_0\M))
\q(\a\ge 0).\leqno(2.4.2)$$
So (2.4.1) follows from (1.4.3).
This finishes the proof of Theorem~1.

\ms\nin
{\bf 2.5.~Proof of Theorem~2.}
In the notation of (1.1) we have
$$x_i=t^{c_{k,i}}x_i^{(k)}.$$
Then $(I_k)_0\subset\V_0$ is given by
$$x_i=0\q\h{for}\,\,\,c_{k,i}>0.$$
However, this is independent of $k\in[1,a]$ by definition
(i.e.\ $c_{k,i}=\min(a_i,k)$).
So the first assertion follows.

For the second assertion we first show that the two spaces have the
same dimension.
With the notation of (1.4) we have
$$\aligned &\Db_{\C}\bigl(\Im(N:F^{j+1}H_{\infty,\C,1}\to
F^jH_{\infty,\C,1})\bigr)\\
&=\Coim(N:H_{\infty,\C,1}/F^{-j}H_{\infty,\C,1}\to H_{\infty,\C,1}/
F^{-j-1}H_{\infty,\C,1})\\
&=H_{\infty,\C,1}/(F^{-j}H_{\infty,\C,1}+\Ker\,N).\endaligned
\leqno(2.5.1)$$
Here the first isomorphism is induced by a polarization, and the
last isomorphism follows from the strict compatibility of
$$N:(H_{\infty,\C,1},F)\to(H_{\infty,\C,1},F[-1]),$$
since $N$ is a morphism of mixed Hodge structure of type
$(-1,-1)$, see [D2].
The last term of (2.5.1) for $j=0$ is further isomorphic to
$$\aligned &\q\frac{H_{\infty,\C,1}/F^0H_{\infty,\C,1}}
{(F^0H_{\infty,\C,1}+\Ker\,N)/F^0H_{\infty,\C,1}}\\
&=\frac{\Gr_V^0\E}{\Im(H^{\rm inv}_{\infty,\C}\to H_{\infty,\C,1}/
F^0H_{\infty,\C,1})}.\endaligned$$
By (2.4.2) and (1.4.3) for $\a=0$ and $j=p=0$, the dimension of the
first term of (2.5.1) for $j=0$ coincides with that of
$$\Gr_V^0(\E/\Db(F_0\M))=\frac{\Gr_V^0\E}{\Gr_V^0\Db(F_0\M)}.$$
The dimension of $\Gr_V^0\Db(F_0\M)$ coincides with the dimension
of the image of $(I_1)_0$ in the unipotent monodromy part since it
can be defined by using the $V$-filtration as in (1.4.2).
So we get the coincidence of the dimensions of the two spaces in
the second assertion, and it is enough to show an inclusion.

The unipotent monodromy part of $(I_1)_0$ is identified with the
image of $\Db(F_0\M)$ in $\Gr_V^0\E/H^{\rm inv}_{\infty,\Z}$.
Here we have to divide $\Gr_V^0\E$ by $\G_0=H^{\rm inv}_{\infty,\Z}$.
To show an inclusion, it is then sufficient to show that the image
of $\Db(F_0\M)$ in $\Gr_V^0\E$ contains the image of
$\Ker\, N\subset H_{\infty,\C,1}$ in
$$\Gr_V^0\E=H_{\infty,\C,1}/F^0H_{\infty,\C,1},$$
i.e.
$$\langle u,v\rangle\in\Oc_{\D}\q\h{if}\,\,\,u\in\Ker\,t\d_t,\,
v\in F_0\M,$$
since $N$ corresponds to $-t\d_t$ on $\Gr_V^0\M$.
(Note that the $V$-filtration in the one-variable case splits
by the action of $t\d_t$.)
But the above assertion follows from (1.3.1).
So the second assertion is proved.
The last assertion then follows from the definition of
the N\'eron model of Green, Griffiths and Kerr [GKK].
This finishes the proof of Theorem~2.

\ms\nin
{\bf 2.6.~Proof of Theorem~3.}
The assertion follows from Theorems~1 and 2 together with
(1.5) and (2.3) by reducing to [Sa2], [Sa3].

\ms\nin
{\bf Remark~2.7.}
In the notation of the last part of the introduction, there is an
injection $J_{\D}^{\rm BPS}(\Hb)\hookrightarrow J_{\D}^{\rm Sch}(\Hb)$
using horizontal sections passing through the monodromy invariant
part of $J_{\D}^C(\Hb)_0$, see [Sch].
Even in the unipotent monodromy case, however, this cannot be
continuous unless it is an isomorphism (e.g. the level is 1).
Here we cannot use admissible normal functions as in the
abelian scheme case explained in the introduction since
the Griffiths transversality gives a strong restriction.
However, the topology of $J_{\D}^{\rm BPS}(\Hb)$ is induced by
the inclusion $J_{\D}^{\rm BPS}(\Hb)\hookrightarrow J_{\D}^C(\Hb)$
by definition [BPS] (in the unipotent monodromy case), and this
implies a contradiction by Theorem~2 if the morphism is continuous
and the surjection $J_{\D}^{\rm Sch}(\Hb)\to J_{\D}^{\rm BPS}(\Hb)
=J_{\D}^{\rm GGK}(\Hb)$ is not bijective.

\ms\nin
{\bf Remark~2.8.}
Theorem~3 shows that the first blow-up $Y_1$ is already a Hausdorff
complex manifold.
As is remarked by the referee, Proposition~2.2 holds for that space
as well, and so $Y_1$, instead of $J_{\D}^{\Sch}(\Hb)^0$, could be
used as the identity component of a N\'eron model.
However, it is not easy to generalize this to the case $\dim S>1$
even in the normal crossing divisor case.

\section{Remarks about the higher dimensional case}

\nin
In this section we give some remarks for the case where the base
space $S$ is not a curve. 

\ms\nin
{\bf Remark~3.1.}
In the case $\dim S>1$, $J_S^{\Sch}(\Hb)$ may have singularities,
caused by the fact that $F_0\M$ is not always locally free.
This can happen even when $\Hb$ is a nilpotent orbit on
$S^*=(\D^*)^2$.
For example, consider the case where the limit mixed Hodge
structure $H_{\infty}$ has rank 4 with type
$$(1,-1),\,\,(-1,1),\,\,(0,-2),\,\,(-2,0),$$
and $N_1,N_2$ are nonzero.
Then $J_S^{\Sch}(\Hb)$ is locally defined by
$$x_1t_1=x_2t_2,$$
in an open neighborhood of $0\in \C^5$ with coordinates
$x_1,x_2,x_3,t_1,t_2$.
Indeed, we have
$$F_0\M=\I_0\oplus\Oc_S,$$
as an $\Oc_S$-module where $\I_0$ denotes the sheaf of ideals of
$0\in S=\D^2$.
(More precisely, $x_1,x_2$ respectively correspond to the two
generators $t_1,t_2$ of $\I_0$, and $x_3$ to the generator $1$ of
$\Oc_S$.)
Hence $F_0\M$ is not a locally free sheaf in this case; in fact,
it is even non-reflexive since
$$\I_0^{\vee}:=\Hom_{\Oc_S}(\I_0,\Oc_S)=\Oc_S,$$
using the short exact sequence $0\to\I_0\to\Oc_S\to \C_0\to 0$.
This calculation implies that the reflexive hull, i.e.\ the double
dual $(F_0\M)^{\vee\vee}$, is free in this case.

\ms\nin
{\bf Remark~3.2.}
The torsion-freeness of $F_0\M$ is not stable by the pull-back under
morphisms of base spaces.
For example, if $\sigma:S'\to S$ is the blow-up along the origin
with $S$, $(\M,F)$ as above, then the pull-back $\sigma^*F_0\M$
has torsion.
Indeed, $\I_0$ is quasi-isomorphic to the mapping cone of
$$(t_1,t_2):\Oc_S\to\Oc_S\oplus\Oc_S,$$
and its pull-back by $\sigma$ is locally the mapping cone of
$$(t'_1,t'_1t'_2):\Oc_{S'}\to\Oc_{S'}\oplus\Oc_{S'},$$
where $t'_1,t'_2$ are local coordinates of $S'$ such that
$\sigma^*t_1=t'_1$, $\sigma^*t_2=t'_1t'_2$.
Then the cokernel of the morphism has nontrivial $t'_1$-torsion.

\bs\nin
{\bf Remark~3.3.}
The freeness of the reflexive hull in the above example
does not hold in general even in the normal crossing case,
e.g.\ if $\Hb$ is a nilpotent orbit of three variables such that
$H_{\infty}$ has dimension 8 with the same type as above and the
images of $F^1H_{\infty}$ by $N_1$, $N_2$, $N_3$ are 1-dimensional
and are not compatible subspaces.
The last condition is equivalent to the condition that they
are distinct to each other and span a 2-dimensional subspace.
So there are $u_1$, $u_2$ such that $u_1$, $u_2$ and $u_1+u_2$
respectively generate the images of $F^1H_{\infty}$ by $N_1$,
$N_2$ and $N_3$.
Then $F_0\M$ is a direct sum of a free $\Oc_S$-module of rank 2
and a coherent $\Oc_S$-module $\M'$ which is generated by
$t_1^{-1}\u_1$, $t_2^{-1}\u_2$ and $t_3^{-1}(\u_1+\u_2)$ where
$S=\D^3$, and $\u_1$, $\u_2$ are defined as in (1.5.1).
This implies that $\M'$ has a a free resolution defined by
$$(t_1,t_2,t_3):\Oc_S\to \Oc_S\oplus\Oc_S\oplus\Oc_S.$$
Indeed, if $\M''$ denote the cokernel of this morphism, then
$H^k_{\{0\}}\M''=0$ for $k\le 1$, and hence $\M''=j_{0*}j_0^*\M''$
where $j_0:S\setminus\{0\}\to S$ is the inclusion.
It is easy to show that $j_0^*\M'=j_0^*\M''$.
So there are morphisms $\M''\to\M'\to j_{0*}j_0^*\M''=\M''$ whose
composition is the identity.
Then the kernel of the surjection $\M'\to\M''$ vanishes, and
we get the isomorphism $\M''=\M'$.

The above resolution is the first two terms of the Koszul complex
$K^{\sb}(\Oc_S;t_1,t_2,t_3)$ which will be denoted by $K^{\sb}$.
This implies that $\M'$ is self-dual, i.e.
$$\Db\M'=\Hc^{-1}\Db(\sigma_{\le 1}K^{\sb})=
\Hc^2(\sigma_{\ge 2}K^{\sb})=\Hc^1(\sigma_{\le 1}K^{\sb})=\M',$$
using the self-duality and the exactness at the middle terms of the
Koszul complex where $\Db(*):=\Hom_{\Oc_S}(*,\Oc_S)$.
(For the truncation $\sigma_{\le p}$, see [D2].)
So $F_0\M$ is self-dual and hence reflexive.
Note that $F_0\M$ cannot be free since the freeness implies a free
resolution of $\C$ over $\Oc_S$ with three terms (shrinking $S$ if
necessary).


\begin{thebibliography}{GGK2}

\bibitem[BBD]{BBD}
A.~Beilinson, J.~Bernstein and P.~Deligne, Faisceaux pervers,
Ast\'erisque, vol. 100, Soc.\ Math.\ France, Paris, 1982.

\bibitem[BPS]{BPS}
P.~Brosnan, G.~Pearlstein and M.~Saito,
A generalization of the N\'eron models of Green, Griffiths and Kerr,
preprint (arXiv:0809.5185).

\bibitem[Cl]{Cl}
H.~Clemens, The N\'eron model for families of intermediate Jacobians
acquiring ``algebraic'' singularities, Publ.\ Math.\ IHES 58 (1983),
5--18.

\bibitem[D1]{D1}
P.~Deligne, Equations diff\'erentielles\`a points singuliers
r\'eguliers, Lect.\ Notes in Math. vol. 163, Springer, Berlin, 1970.

\bibitem[D2]{D2}
P.~Deligne, Th\'eorie de Hodge II, Publ. Math. IHES 40 (1971),
5--58.

\bibitem[EZ]{EZ}
F.~El Zein and S.~Zucker, Extendability of normal functions
associated to algebraic cycles, in Topics in transcendental
algebraic geometry, Ann. Math. Stud., 106, Princeton Univ. Press,
Princeton, N.J., 1984, pp.~269--288.

\bibitem[GGK1]{}
M.~Green, P.~Griffiths and M.~Kerr,
Neron models and boundary components for degenerations of Hodge
structures of mirror quintic type, 
in ``Curves and Abelian Varieties (V. Alexeev, Ed.)'',
Contemp.\ Math.\ 465 (2007), AMS, 71-145.

\bibitem[GGK2]{GGK2}
M.~Green, P.~Griffiths and M.~Kerr,
N\'eron models and limits of Abel-Jacobi mappings,
Compositio Math. 146 (2010), 288--366.

\bibitem[OSS]{OSS}
C.~Okonek, M.~Schneider and H.~Spindler,
Vector bundles on complex projective spaces, Birkh\"auser, 1980.

\bibitem[Sa1]{Sa1}
M.~Saito, Modules de Hodge polarisables, Publ. RIMS, Kyoto Univ.
24 (1988), 849--995.

\bibitem[Sa2]{Sa2}
M.~Saito, Admissible normal functions, J.\ Algebraic Geom.\ 5 (1996),
235--276.

\bibitem[Sa3]{Sa3}
M.~Saito, Hausdorff property of the N\'eron models of Green,
Griffiths and Kerr, preprint (arXiv:0803.2771).

\bibitem[Sch]{Sch}
C.~Schnell, Complex analytic N\'eron models for arbitrary families
of intermediate Jacobians, preprint (arXiv:0910.0662).

\bibitem[Yo]{Yo}
A.~Young, Complex analytic N\'eron models for degenerating Abelian
varieties over higher dimensional parameter spaces,
Ph.\ D.\ Thesis, Princeton University, Sept.\ 2008.

\bibitem[Zu]{Zu}
S.~Zucker, Generalized intermediate Jacobians and the theorem on
normal functions, Inv.\ Math.\ 33 (1976),185--222.

\end{thebibliography}
\end{document}